\documentclass[12pt,oneside,reqno]{amsart}
\usepackage{txfonts}
\usepackage{bbm}
\usepackage{amsmath}
\usepackage{graphicx}
\usepackage{mathrsfs}
\usepackage{stmaryrd}
\usepackage{amsfonts}
\usepackage{enumerate,amsmath,amssymb,amsthm}

\numberwithin{equation}{section}

\newcommand{\be}{\begin{eqnarray}}
\newcommand{\ee}{\end{eqnarray}}
\newcommand{\ce}{\begin{eqnarray*}}
\newcommand{\de}{\end{eqnarray*}}
\newtheorem{theorem}{Theorem}[section]
\newtheorem{lemma}[theorem]{Lemma}
\newtheorem{remark}[theorem]{Remark}
\newtheorem{definition}[theorem]{Definition}
\newtheorem{proposition}[theorem]{Proposition}
\newtheorem{Examples}[theorem]{Example}
\newtheorem{corollary}[theorem]{Corollary}

\def\eps{\varepsilon}

\def\p{\partial}

\def\[{{\Big[}}
\def\]{{\Big]}}
\def\<{{\langle}}
\def\>{{\rangle}}
\def\({{\Big(}}
\def\){{\Big)}}

\def\bx{{\mathbf{x}}}

\def\dif{{\mathord{{\rm d}}}}

\def\no{\nonumber}
\def\={&\!\!=\!\!&}
\def\bt{\begin{theorem}}
\def\et{\end{theorem}}
\def\bl{\begin{lemma}}
\def\el{\end{lemma}}
\def\br{\begin{remark}}
\def\er{\end{remark}}

\def\bd{\begin{definition}}
\def\ed{\end{definition}}
\def\bp{\begin{proposition}}
\def\ep{\end{proposition}}
\def\bc{\begin{corollary}}
\def\ec{\end{corollary}}
\def\bx{\begin{Examples}}
\def\ex{\end{Examples}}

\def\cM{{\mathcal M}}

\def\cP{{\mathcal P}}

\def\mB{{\mathbb B}}

\def\mE{{\mathbb E}}

\def\mH{{\mathbb H}}

\def\mM{{\mathbb M}}
\def\mN{{\mathbb N}}

\def\mP{{\mathbb P}}

\def\mR{{\mathbb R}}

\def\mT{{\mathbb T}}
\def\mU{{\mathbb U}}

\def\mX{{\mathbb X}}

\def\sF{{\mathscr F}}

\def\geq{\geqslant}
\def\leq{\leqslant}

\allowdisplaybreaks

\begin{document}

\title{Exponential Ergodicity of stochastic Burgers equations driven by $\alpha$-stable processes}

\date{}
\author{Zhao Dong, Lihu Xu and Xicheng Zhang}

\address{Zhao Dong: Institute of Applied Mathematics,
Academy of Mathematics and Systems Sciences, Academia Sinica, P.R.China\\
Email: dzhao@amt.ac.cn}

\address{Lihu Xu: Department of Mathematics, Brunel University,
Kingston Lane, Uxbridge, Middlesex UB8 3PH, United Kingdom, Email: Lihu.Xu@brunel.ac.uk}

\address{Xicheng Zhang:
School of Mathematics and Statistics, Wuhan University,
Wuhan, Hubei 430072, P.R.China\\
Email: XichengZhang@gmail.com
 }

\begin{abstract}
In this work, we prove the strong Feller property and the exponential ergodicity of stochastic Burgers equations driven
by $\alpha/2$-subordinated cylindrical Brownian motions with $\alpha\in(1,2)$.
To prove the results, we truncate the nonlinearity and
use the derivative formula for SDEs driven by $\alpha$-stable noises established in \cite{Zh}.
\end{abstract}

\maketitle
\rm

\section{Introduction}

Stochastic Burgers and Navier-Stokes equations, as models of studying the statistic theory of the turbulent fluid motion,
has been studied in many literatures in past twenty years. In particular, the existence-uniqueness and ergodicity have
been studied by many authors under non-degenerate or degenerate random perturbations
(cf. \cite{Be-Ca-Jo,Da-De-Te,Gy-Nu,Fl-Ma,Go-Ma,Ha-Ma,Ro-Xu} etc.). In these works, the random forces
are assumed to be the Brownian noise, which can be naturally regarded as a continuous time model.

In recent years, the stochastic equations driven by L\'evy type noises also attract much attention
(cf. \cite{Do-Xu}, \cite{Do-Xi}, \cite{Pe-Za}-\cite{Pr-Za}, \cite{Wa}-\cite{Zh}, etc.).
It was proved in \cite{Do-Xu} and \cite{Do-Xi} that there is a unique invariant measure
for stochastic Burgers and 2D Navier-Stokes equations with L\'evy noises.
In these two works, the L\'evy noises are assumed to be square integrable.
This restriction clearly rules out the interesting $\alpha$-stable noises. It should be stressed that
since the $\alpha$-stable noise exhibits the heavy tailed phenomenon, the stochastic equation
driven by $\alpha$-stable processes recently causes great interest in physics
(cf. \cite{Ch-Go-Kl-Me-Ta,Lu-La,Ma-St,So-We-Sw} etc.).

We shall consider in this paper the following stochastic Burgers equation on torus $\mathbb{T}=\mathbb{R}/(2\pi\mathbb{Z})$:
\begin{align}
\p_t u_t=u''_t-u_tu_t'+\dot\xi_t,\label{NSE}
\end{align}
where $\dot\xi_t$ is some time-white noise. As mentioned above, when $\dot\xi_t$ is additive Brownian noise,
this type of equation has been intensively studied. In \cite{Be-Ca-Jo},
Bertini, Cancrini and Jona-Lasini used Cole-Hopf's transformation to reduce equation (\ref{NSE}) to a
linear heat equation and obtained the existence of solutions. In \cite{Da-Za}, the ergodicity was also proved by using some
truncation technique (see also \cite{Fl-Ma,Go-Ma,Ha-Ma} etc. for stochastic Navier-Stokes equations).
In the present work we shall assume that $\dot\xi_t$ is a type of $\alpha$-stable noise
called $\alpha/2$-subordinated cylindrical Brownian noise
and prove the exponential ergodicity of equation \eqref{NSE}. There have been some
results on ergodicity of stochastic systems driven by
$\alpha$-stable type noises (cf.~\cite{Xu-Ze1,Pr-Sh-Xu-Za, Ku, Xu}).
In \cite{Ku}, Kulik obtained a nice
criterion for the exponential mixing of a family of SDEs driven by $\alpha$-stable noises.
 We refer to \cite{Xu-Ze1} for the exponential mixing of stochastic spin systems with
$\alpha$-stable noises, and to \cite{Pr-Sh-Xu-Za} for the exponential mixing
of a family of semi-linear SPDEs with Lipschitz nonlinearity.

Let us now discuss the approach to the ergodicity. In a previous work \cite{Do-Xu-Zh},
we have proved the existence of invariant measures for stochastic 2D Navier-Stokes equation by
estimating the fractional moments. The proofs clearly also works for Burgers equation (\ref{NSE}).
To prove the exponential ergodicity, we shall use the Harris theorem (cf. \cite{Ha-Ma-Sh}).
Thus, the main task is to
verify the conditions in Harris theorem, where an important step in our proof
is to prove the strong Feller property
for truncated equation. It is well known that the truncating nonlinearity technique is a usual tool to establish the strong Feller
property for Navier-Stokes and Ginzburg-Landau type equations (\cite{Fl-Ma, Ec-Ha, Ro-Xu, Xu}).
To prove the strong Feller property, we shall truncate the quadratic nonlinearity of equation \eqref{NSE} and
apply a derivative formula established in \cite{Zh}.

This paper is organized as follows: In Section 2, we give some necessary notions and notations. In particular,
we study the stochastic convolutions in Hilbert space about the $\alpha/2$-subordinated cylindrical Brownian motions.
In Section 3, we present a general result about the strong Feller property for SPDEs driven by
$\alpha/2$-subordinated cylindrical Brownian motions. This result generalizes the corresponding one
in \cite[Theorem 4.1]{Zh}. In Section 4, we prove our main result Theorem \ref{Main} by using
suitable truncation technique and verifying the Harris conditions. In appendix, we study a deterministic Burgers equation and give some necessary dependence relation
about the initial values. The result is by no means new. Since the proof is not so long,
we include it here for the reader's convenience.

We conclude this section by introducing the following conventions: The letter $C$ with or without subscripts
will denote an unimportant constant, whose value may change in different occasions. Moreover,
let $\mathbb U$ be a Banach space, for $R>0$ we shall denote the ball in $\mU$ by
$$
\mB^\mU_R:=\{u\in\mU: \|u\|_\mU\leq R\}.
$$


\section{Preliminaries}

Let $\mH$ be a real separable Hilbert space with the inner product $\<\cdot,\cdot\>_0$. The norm in $\mH$
is denoted by $\|\cdot\|_0$. Let $A$ be a positive self-adjoint operator on $\mH$ with discrete spectral, i.e.,
there exists an orthogonal basis $\{e_k\}_{k\in\mN}$ and a sequence of real numbers
$0<\lambda_1\leq\lambda_2\leq\cdots\leq\lambda_k\to\infty$ such that
$$
A e_k=\lambda_k e_k.
$$
For $\gamma\in\mR$, let $\mH^\gamma$ be the domain of the fractional operator $A^{\frac{\gamma}{2}}$, i.e.,
$$
\mH^\gamma:=A^{-\frac{\gamma}{2}}(\mH)
=\left\{\sum_k\lambda_k^{-\frac{\gamma}{2}}a_ke_k: (a_k)_{k\in\mN}\subset\mR, \sum_k a_k^2<+\infty\right\},
$$
with the inner product
$$
\<u,v\>_\gamma:=\<A^{\frac{\gamma}{2}}u,A^{\frac{\gamma}{2}}v\>_0=\sum_k\lambda_k^\gamma\<u,e_k\>_0\<v,e_k\>_0.
$$
The semigroup associated to $A$ is defined by
$$
\mathrm{e}^{-tA}u:=\sum_k \mathrm{e}^{-t\lambda_k}\<u,e_k\>_0 e_k,\ \ t\geq 0.
$$
It is easy to see that for any $\gamma>0$,
\begin{align}
\|A^\gamma \mathrm{e}^{-tA}u\|_0\leq \sup_{x>0}(x^\gamma
\mathrm{e}^{-x})t^{-\gamma}\|u\|_0=\gamma^\gamma \mathrm{e}^{-\gamma}t^{-\gamma}\|u\|_0.\label{Semi}
\end{align}
For a sequence of bounded real numbers
$\beta=(\beta_k)_{k\in\mN}$, let us define
$$
Q_\beta: \mH\to\mH; \ \ Q_\beta u:=\sum_{k=1}^\infty\beta_k \<u,e_k\>_0 e_k.
$$
\bl
Suppose that for some $\delta>0$ and $\theta,\theta'\in\mR$ with $\theta>\theta'$,
\begin{align}
\delta\lambda^{-\frac{\theta}{2}}_k\leq|\beta_k|\leq \delta^{-1}\lambda^{-\frac{\theta'}{2}}_k,\ \ \forall k\in\mN,\label{EW1}
\end{align}
Then we have
\begin{align}
\|A^{\frac{\theta'}{2}}Q_\beta u\|_0\leq \delta^{-1}\|u\|_0,\ \ u\in\mH^0\label{EE6}
\end{align}
and
\begin{align}
\|Q^{-1}_\beta u\|_0\leq \delta^{-1}\|A^{\frac{\theta}{2}}u\|_0,\ \ u\in\mH^\theta.\label{EE7}
\end{align}
\el
\begin{proof}
By definition, we have
$$
\|A^{\frac{\theta'}{2}}Q_\beta u\|^2_0=\sum_k|\beta_k|^2\lambda_k^{\theta'}\<u,e_k\>^2_0
\leq \delta^{-2}\sum_k\<u,e_k\>^2_0=\delta^{-2}\|u\|^2_0,
$$
and
$$
\|Q^{-1}_\beta u\|^2_0=\sum_{k=1}^\infty|\beta_k|^{-2}\<u,e_k\>^2_0
\leq \delta^{-2}\sum_k\lambda_k^{\theta}\<u,e_k\>^2_0
=\delta^{-2}\|A^{\frac{\theta}{2}}u\|^2_0.
$$
The estimates follow.
\end{proof}

Let $\{W^k_t, t\geq 0\}_{k\in\mN}$ be a sequence of independent standard one-dimensional Brownian motion
on some probability space $(\Omega,\sF,\mP)$.
The cylindrical Brownian motion on $\mH$ is defined by
$$
W_t:=\sum_k W^k_te_k.
$$
For $\alpha\in(0,2)$, let $S_t$ be an independent $\alpha/2$-stable subordinator, i.e., an increasing one
dimensional L\'evy process with Laplace transform
$$
\mE \mathrm{e}^{-\eta S_t}=\mathrm{e}^{-t|\eta|^{\alpha/2}},\ \ \eta>0.
$$
The subordinated cylindrical Brownian motion $\{L_t\}_{t\geq 0}$ on $\mH$ is defined by
$$
L_t:=W_{S_t}.
$$
Notice that in general $L_t$ does not belong to $\mH$.

We recall the following estimate about the subordinator $S_t$.
\bl
We have
\begin{align}
\mP(S_t\leq r)>0,\ \ r,t>0,\label{EE9}
\end{align}
and
\begin{align}
\mE\Big(S_t^{-q}\Big)\leq C t^{-\frac{2q}{\alpha}},\ \ q,t>0.\label{EE10}
\end{align}
\el
\begin{proof}
Estimate (\ref{EE9}) follows by the strict positivity of the distributional density $p_t(s)$ of $S_t$.
For (\ref{EE10}), recalling that $p_t(s)$ satisfies (cf. \cite[(14)]{Bo-St-Sz})
$$
p_t(s)\leq C ts^{-1-\frac{\alpha}{2}}\mathrm{e}^{-ts^{-\frac{\alpha}{2}}},
$$
we have
\begin{align*}
\mE\left(S^{-q}_t\right)
\leq C\int^\infty_0ts^{-1-\frac{\alpha+2q}{2}}\mathrm{e}^{-ts^{-\frac{\alpha}{2}}}\dif s
=C t^{-\frac{2q}{\alpha}}\int^\infty_0 u^{\frac{2q}{\alpha}}\mathrm{e}^{-u}\dif u,
\end{align*}
where the last equality is due to the change of variable $u=ts^{-\frac{\alpha}{2}}$, and
$C$ only depends on $\alpha,q$.
\end{proof}

Let us now consider the following stochastic convolution:
$$
Z_t:=\int^t_0\mathrm{e}^{-(t-s)A}Q_\beta\dif L_s=\sum_k\int^t_0 \mathrm{e}^{-(t-s)\lambda_k}\beta_k\dif W^k_{S_s} e_k,
$$
where $Q_\beta$ denotes the intensity of the noise.
The following estimate about $Z_t$ will play an important role in the next sections (cf. \cite{Pr-Za, Pe-Za}).
\bl
Suppose that for some $\gamma\in\mR$,
\begin{align}
K_\gamma:=\sum_k\lambda_k^\gamma|\beta_k|^2<+\infty.\label{DD}
\end{align}
Then for any $p\in(0,\alpha)$ and $T>0$,
\begin{align}
\sup_{t\in[0,T]}\mE\|Z_t\|_{\gamma+1}^p\leq C_{\alpha,p} K_\gamma^{\frac{p}{2}}T^{\frac{p}{\alpha}-\frac{p}{2}},\label{Es1}
\end{align}
and for any $\theta<\gamma$,
\begin{align}
\mE\left(\sup_{t\in[0,T]}\|Z_t\|_\theta^p\right)
\leq C_{\alpha,p}K_\gamma^{\frac{p}{2}}T^{\frac{p}{\alpha}}\left(1+T^{\frac{\gamma-\theta}{2}}\right),\label{Ep1}
\end{align}
and for any $\eps>0$,
\begin{align}
\mP\left(\sup_{t\in[0,T]}\|Z_t\|_\theta\leq\eps\right)>0.\label{Ep11}
\end{align}
Moreover, $t\mapsto Z_t$ is almost surely c\`adl\`ag in $\mH^\theta$.
\el
\begin{proof}
Estimate (\ref{Es1}) follows by \cite[Proposition 4.2]{Zh}. Next, we prove (\ref{Ep1}).
For any $p\in(0,\alpha)$, by Burkh\"older's inequality for Brownian motion, we have
\begin{align}
\mE\left(\sup_{t\in[0,T]}\|A^{\frac{\gamma}{2}} Q_\beta L_t\|_0^p\right)
&=\mE\left(\mE\left(\sup_{t\in[0,T]}\|A^{\frac{\gamma}{2}} Q_\beta W_{\ell_t}\|_0^p\right)\Bigg|_{\ell=S}\right)\no\\
&\leq\mE\left(\mE\left(\sup_{s\in[0,\ell_T]}\|A^{\frac{\gamma}{2}} Q_\beta W_s\|_0^p\right)\Bigg|_{\ell=S}\right)\no\\
&\leq C_p\mE\left(\|A^{\frac{\gamma}{2}} Q_\beta\|_{\mathrm{H.S.}}^pS_T^{\frac{p}{2}}\right)\no\\
&=C_p\left(\sum_k\lambda^\gamma_k\beta_k^2\right)^{\frac{p}{2}}\mE\left(S_1^{\frac{p}{2}}\right)T^{\frac{p}{\alpha}}.\label{Es2}
\end{align}
In particular,
\begin{align}
\mbox{$t\mapsto Q_\beta L_t$ is almost surely c\`adl\`ag in $\mH^\gamma$}.\label{ET3}
\end{align}
On the other hand, by integration by parts formula, we have
\begin{align}
Z_t=Q_\beta L_t+\int^t_0A\mathrm{e}^{-(t-s)A}Q_\beta L_s\dif s.\label{EW23}
\end{align}
Hence, for any $\theta<\gamma$, by (\ref{Semi}) we have
\begin{align}
\|A^{\frac{\theta}{2}}Z_t\|_0&\leq \|A^{\frac{\theta}{2}} Q_\beta L_t\|_0
+\int^t_0\|A^{1+\frac{\theta-\gamma}{2}}\mathrm{e}^{-(t-s)A}A^{\frac{\gamma}{2}}Q_\beta L_s\|_0\dif s\no\\
&\leq \lambda_1^{\theta-\gamma}\|A^{\frac{\gamma}{2}} Q_\beta L_t\|_0+C\int^t_0\frac{\|A^{\frac{\gamma}{2}}Q_\beta L_s\|_0}{(t-s)^{1+\frac{\theta-\gamma}{2}}}\dif s\no\\
&\leq C\sup_{s\in[0,t]}\|A^{\frac{\gamma}{2}}Q_\beta L_s\|_0\left(1+t^{\frac{\gamma-\theta}{2}}\right)
=:\sup_{s\in[0,t]}\|A^{\frac{\gamma}{2}}Q_\beta L_s\|_0\cdot \eta_t.\label{Es3}
\end{align}
Estimate (\ref{Ep1}) then follows by combining (\ref{Es2}) and (\ref{Es3}). Moreover, we also have
that $t\mapsto\int^t_0A\mathrm{e}^{-(t-s)A}Q_\beta L_s\dif s$ is continuous in $\mH^\theta$.
Thus, the c\`adl\`ag property of $t\mapsto Z_t$ in $\mH^\theta$ follows by (\ref{ET3}) and (\ref{EW23}).

Now, we prove (\ref{Ep11}). By (\ref{Es3}) we have
\begin{align*}
\mP\left(\sup_{t\in[0,T]}\|Z_t\|_\theta\leq\eps\right)&\geq
\mP\left(\sup_{t\in[0,T]}\|A^{\frac{\gamma}{2}}Q_\beta W_{S_t}\|_0\leq\eps\eta_T^{-1}\right)\\
&\geq\mP\left(\sup_{t\in[0,S_T]}\|A^{\frac{\gamma}{2}}Q_\beta W_t\|_0\leq\eps\eta_T^{-1}\right)\\
&\geq\mP\left(\sup_{t\in[0,S_T]}\|A^{\frac{\gamma}{2}}Q_\beta W_t\|_0\leq\eps\eta_T^{-1}; S_T\leq 1\right)\\
&\geq\mP\left(\sup_{t\in[0,1]}\|A^{\frac{\gamma}{2}}Q_\beta W_t\|_0\leq\eps\eta_T^{-1}; S_T\leq 1\right)\\
&=\mP\left(\sup_{t\in[0,1]}\|A^{\frac{\gamma}{2}}Q_\beta W_t\|_0\leq\eps\eta_T^{-1}\right)\mP(S_T\leq 1)>0.
\end{align*}
The last step is due to the fact that each term is positive.
\end{proof}

\section{Strong Feller property of SPDEs driven by subordinated cylindrical Brownian motions}

In this section, we consider the following general SPDE in Hilbert space $\mH$:
\begin{align}
\dif u_t=[-Au_t+F(u_t)]\dif t+Q_\beta\dif L_t,\ \ u_0=\varphi\in\mH,\label{Eq}
\end{align}
where for some $\delta>0$ and $\theta\geq \theta'\geq 0$,
\begin{align}
\delta\lambda^{-\frac{\theta}{2}}_k\leq|\beta_k|\leq \delta^{-1}\lambda^{-\frac{\theta'}{2}}_k,\ \
\forall k\in\mN,\label{EW11}
\end{align}
and for some $\gamma,\gamma'\geq 0$,
\begin{align}
F:\mH^\gamma\to\mH^{-\gamma'}\mbox{ is bounded and Lipschitz continuous}.\label{EW22}
\end{align}

We need the following important constant:
\begin{align}
\theta_0:=\inf\left\{\theta>0: \sum_k\lambda_k^{-\theta}<+\infty\right\}.\label{ET2}
\end{align}
The aim of this section is to prove that
\bt\label{Th}
Let $\alpha\in(1,2)$ and $Z_t:=\int^t_0 \mathrm{e}^{-(t-s)A}Q_\beta\dif L_s$.
Assume that (\ref{EW11}) and (\ref{EW22}) hold with
\begin{align}
\gamma-\theta'<1-\theta_0,\ \gamma+\gamma'<2,\label{Es8}
\end{align}
then for any $\varphi\in\mH$, there exists a unique $u_t=u_t(\varphi)$ satisfying that
$$
u_t-Z_t\in C([0,\infty);\mH)\cap C((0,\infty);\mH^\gamma),
$$
and
\begin{align}
u_t=\mathrm{e}^{-tA}\varphi+\int^t_0 \mathrm{e}^{-(t-s)A}F(u_s)\dif s+Z_t.\label{Es4}
\end{align}
If in addition that for some $\sigma\geq 0$,
$$
\gamma\leq\theta<\sigma+\tfrac{2}{\alpha},\ \ \theta+\gamma'<2,
$$
then for any bounded Borel measurable function $\Phi: \mH\to\mR$, $\varphi_1,\varphi_2\in\mH^\sigma$ and $t>0$,
\begin{align}
|\mE \Phi(u_t(\varphi_1))-\mE \Phi(u_t(\varphi_2))|\leq C_t
t^{-\frac{1}{\alpha}-\frac{\theta-\sigma}{2}}\|\Phi\|_\infty\|\varphi_1-\varphi_2\|_\sigma,\label{For5}
\end{align}
where $t\mapsto C_t$ is a continuous increasing function on $[0,\infty)$.
\et

\begin{proof}
The proof is divided into four steps.

(Step 1). We first establish the existence and uniqueness for (\ref{Es4}).
Set $w_t:=u_t-Z_t$. Thus, to solve equation (\ref{Es4}), it suffices to solve the following deterministic equation:
$$
w_t=\mathrm{e}^{-tA}\varphi+\int^t_0\mathrm{e}^{-(t-s)A}F(w_s+Z_s)\dif s.
$$
By (\ref{Es1}), (\ref{EW11}) and (\ref{Es8}), we have
$$
\int^T_0\mE\|Z_t\|_\gamma^p\dif t\leq C_T K^{\frac{p}{2}}_{\gamma-1}<+\infty,\ \ \forall T>0,
$$
where $K_{\gamma-1}$ is defined by (\ref{DD}).
Therefore, there exists a null set $\Omega_0\subset\Omega$ such that for all $\omega\notin \Omega_0$,
$$
Z_t(\omega)\in\mH^\gamma\mbox{ for Lebesgue almost all $t\geq 0$}.
$$
Below, we fix such an $\omega$ and use the standard Picard's iteration argument
to prove the existence. Define $w^{(0)}_t:=\mathrm{e}^{-tA}\varphi$ and for $n\in\mN$,
\begin{align}
w^{(n)}_t:=\mathrm{e}^{-tA}\varphi+\int^t_0\mathrm{e}^{-(t-s)A}F(w^{(n-1)}_s+Z_s)\dif s.\label{Es5}
\end{align}
By (\ref{Semi}), we have
\begin{align}
\|w^{(n)}_t\|_\gamma&\leq \|A^{\frac{\gamma}{2}}\mathrm{e}^{-tA}\varphi\|_0+\int^t_0\|A^{\frac{\gamma+\gamma'}{2}}
\mathrm{e}^{-(t-s)A}A^{-\frac{\gamma'}{2}}F(w^{(n-1)}_s+Z_s)\|_0\dif s\no\\
&\leq Ct^{-\frac{\gamma}{2}}\|\varphi\|_0+C\int^t_0(t-s)^{-\frac{\gamma+\gamma'}{2}}
\|F(w^{(n-1)}_s+Z_s)\|_{-\gamma'}\dif s\no\\
&\leq Ct^{-\frac{\gamma}{2}}\|\varphi\|_0+C\sup_{u\in\mH^\gamma}\|F(u)\|_{-\gamma'}
\int^t_0(t-s)^{-\frac{\gamma+\gamma'}{2}}\dif s\no\\
&= Ct^{-\frac{\gamma}{2}}\|\varphi\|_0+Ct^{1-\frac{\gamma+\gamma'}{2}}\sup_{u\in\mH^\gamma}\|F(u)\|_{-\gamma'}.\label{EE1}
\end{align}
Similarly, for any $n,m\in\mN$, we also have
\begin{align*}
\|w^{(n)}_t-w^{(m)}_t\|_\gamma&\leq \int^t_0\|A^{\frac{\gamma+\gamma'}{2}}\mathrm{e}^{-(t-s)A}
A^{-\frac{\gamma'}{2}}(F(w^{(n-1)}_s+Z_s)-F(w^{(m-1)}_s+Z_s))\|_0\dif s\\
&\leq C\int^t_0(t-s)^{-\frac{\gamma+\gamma'}{2}}
\|F(w^{(n-1)}_s+Z_s)-F(w^{(m-1)}_s+Z_s)\|_{-\gamma'}\dif s\\
&\leq C\|F\|_{\mathrm{Lip}}\int^t_0(t-s)^{-\frac{\gamma+\gamma'}{2}}\|w^{(n-1)}_s-w^{(m-1)}_s\|_\gamma\dif s,
\end{align*}
where $\|F\|_{\mathrm{Lip}}:=\sup_{u\not=v\in\mH^\gamma}\frac{\|F(u)-F(v)\|_{-\gamma'}}{\|u-v\|_\gamma}$.
This implies that for $q<\tfrac{2}{\gamma+\gamma'}$, $p=\frac{q}{q-1}$ and all $t\in[0,T]$,
\begin{align*}
t^{\frac{\gamma}{2}}\|w^{(n)}_t-w^{(m)}_t\|_\gamma
&\leq C t^{\frac{\gamma}{2}}\left(\int^t_0\left((t-s)^{-\frac{\gamma+\gamma'}{2}}s^{-\frac{\gamma}{2}}\right)^q\dif s\right)^{\frac{1}{q}}
\left(\int^t_0\left(s^{\frac{\gamma}{2}}\|w^{(n-1)}_s-w^{(m-1)}_s\|_\gamma\right)^p\dif s\right)^{\frac{1}{p}}\\
&\leq C t^{\frac{1}{q}-\frac{\gamma+\gamma'}{2}}
\left(\int^t_0\left(s^{\frac{\gamma}{2}}\|w^{(n-1)}_s-w^{(m-1)}_s\|_\gamma\right)^p\dif s\right)^{\frac{1}{p}}.
\end{align*}
Thus, by (\ref{EE1}) and Fatou's lemma, we have
\begin{align*}
\varlimsup_{n,m\to\infty}\sup_{s\in[0,t]}\left(s^{\frac{\gamma}{2}}\|w^{(n)}_s-w^{(m)}_s\|_\gamma\right)^p
\leq C_T\int^t_0\varlimsup_{n,m\to\infty}
\sup_{r\in[0,s]}\left(r^{\frac{\gamma}{2}}\|w^{(n-1)}_r-w^{(m-1)}_r\|_\gamma\right)^p\dif s.
\end{align*}
By Gronwall's inequality, we obtain
\begin{align}
\varlimsup_{n,m\to\infty}\sup_{s\in[0,T]}s^{\frac{\gamma}{2}}\|w^{(n)}_s-w^{(m)}_s\|_\gamma=0.\label{Es7}
\end{align}
Hence, there exists a $w\in C((0,\infty);\mH^\gamma)$ such that for all $T>0$,
$$
\varlimsup_{n\to\infty}\sup_{s\in[0,T]}s^{\frac{\gamma}{2}}\|w^{(n)}_s-w_s\|_\gamma=0.
$$
Taking limits for equation (\ref{Es5}), we obtain the existence of a solution. The uniqueness follows from similar calculations.

(Step 2). Let $\mH_n$ be the finite dimensional subspace of $\mH$ spanned by $\{e_1,\cdots,e_n\}$. Below we always use
the isomorphism:
$$
\mH_n\simeq\mR^n: u=\sum_{k=1}^nu_k e_k,\ \ (u_1,\cdots, u_n)\in\mR^n.
$$
Let $\Pi_n$ be the projection operator from $\mH$ to $\mH_n$
defined by
$$
\Pi_n u:=\sum_{k=1}^n\<u,e_k\>_\mH e_k.
$$
Let $\rho_n$ be a sequence of nonnegative smooth functions with
$$
\mathrm{supp}(\rho_n)\subset\{z\in\mH_n: |z|\leq 1/n\},\ \ \int_{\mH_n}\rho_n(z)\dif z=1.
$$
Define
$$
F_n(u):=\int_{\mH_n}\rho_n(A^{\frac{\gamma}{2}}(u-z))\Pi_n F(z)\dif z=\int_{\mH_n}\rho_n(z)\Pi_n F(u-A^{-\frac{\gamma}{2}}z)\dif z,\ \ u\in\mH_n.
$$
Then
$$
A^{-\frac{\gamma'}{2}}F_n(u)=\int_{\mH_n}\rho_n(u-z)\Pi_nA^{-\frac{\gamma'}{2}} F(z)\dif z.
$$
Since $F:\mH^{\gamma}\to\mH^{-\gamma'}$ is Lipschitz continuous, it is easy to see that
\begin{align}
\sup_{u\in\mH_n}\|\nabla_h A^{-\frac{\gamma'}{2}}F_n(u)\|_0\leq
\sup_{u\not=v}\frac{\|F(u)-F(v)\|_{-\gamma'}}{\|u-v\|_\gamma}\|A^{\frac{\gamma}{2}} h\|_0,\ \ h\in\mH_n.\label{ET1}
\end{align}

Let
$$
L^{(n)}_t:=\sum_{k=1}^nW^k_{S_t}e_k.
$$
Consider the following finite dimensional SDE:
$$
\dif u^{(n)}_t=[-Au^{(n)}_t+F_n(u^{(n)}_t)]\dif t+Q_\beta\dif L^{(n)}_t,\ \ u^{(n)}_0=\varphi\in\mH_n.
$$
By Duhamel's formula, we have
$$
u^{(n)}_t(\varphi)=\mathrm{e}^{-tA}\varphi+\int^t_0 \mathrm{e}^{-(t-s)A}F_n(u^{(n)}_s(\varphi))\dif s
+\int^t_0\mathrm{e}^{-(t-s)A}Q_\beta\dif L^{(n)}_s.
$$
It is easy to see that the directional derivative of $\varphi\mapsto u^{(n)}_t(\varphi)$ along the direction $h\in\mH_n$
satisfies
\begin{align*}
\nabla_hu^{(n)}_t(\varphi)=\mathrm{e}^{-tA}h+\int^t_0 \mathrm{e}^{-(t-s)A}\nabla_h(F_n\circ u^{(n)}_s)(\varphi)\dif s
\end{align*}
By (\ref{ET1}), we further have
\begin{align*}
\|A^{\frac{\theta}{2}}\nabla_hu^{(n)}_t(\varphi)\|_0&\leq \|A^{\frac{\theta}{2}} \mathrm{e}^{-tA}h\|_0
+\int^t_0 \|A^{\frac{\theta+\gamma'}{2}} \mathrm{e}^{-(t-s)A}\nabla_hA^{-\frac{\gamma'}{2}}(F_n\circ u^{(n)}_s)(\varphi)\|_0\dif s\\
&\leq Ct^{\frac{\sigma-\theta}{2}}\|A^{\frac{\sigma}{2}} h\|_0+C\int^t_0 (t-s)^{-\frac{\theta+\gamma'}{2}}
\|\nabla_hA^{-\frac{\gamma'}{2}}(F_n\circ u^{(n)}_s)(\varphi)\|_0\dif s\\
&\leq Ct^{\frac{\sigma-\theta}{2}}\|h\|_\sigma+C\int^t_0 (t-s)^{-\frac{\theta+\gamma'}{2}}
\|A^{\frac{\gamma}{2}}\nabla_h  u^{(n)}_s(\varphi)\|_0\dif s,
\end{align*}
in view of $\gamma\leq\theta$, which implies that
$$
t^{\frac{\theta-\sigma}{2}}\|A^{\frac{\theta}{2}}\nabla_hu^{(n)}_t(\varphi)\|_0\leq C\|h\|_\sigma+
Ct^{\frac{\theta-\sigma}{2}}\int^t_0 \Big((t-s)^{-\frac{\theta+\gamma'}{2}}s^{-\frac{\theta-\sigma}{2}}\Big)
s^{\frac{\theta-\sigma}{2}}\|A^{\frac{\theta}{2}}\nabla_h  u^{(n)}_s(\varphi)\|_0\dif s.
$$
As in the proof of (\ref{Es7}),  we have
\begin{align}
t^{\frac{\theta-\sigma}{2}}\|A^{\frac{\theta}{2}}\nabla_hu^{(n)}_t(\varphi)\|_0\leq C_T\|h\|_\sigma,\ \ h\in\mH_n,\ t\in(0,T],\label{EE8}
\end{align}
where $C_T$ is independent of $n$.

Now, by \cite[Theorem 1.1]{Zh}, we have
$$
\nabla_h\mE \Phi(u^{(n)}_t(\varphi))=\mE\left(\Phi(u^{(n)}_t(\varphi))
\frac{1}{S_t}\int^t_0\<Q^{-1}_\beta\nabla_hu^{(n)}_s(\varphi),\dif L^{(n)}_s\>_0\right).
$$
By H\"older's inequality, (\ref{EE10}) and \cite[Theorem 3.2]{Zh},
for any $p\in(1,\alpha)$ and $q=\frac{p}{p-1}$, we have
\begin{align*}
\|\nabla_h\mE \Phi(u^{(n)}_t(\varphi))\|_0
&\leq\|\Phi\|_\infty\left(\mE\left(\frac{1}{S^{q}_t}\right)\right)^{1/q}
\left(\mE\left|\int^t_0\<Q^{-1}_\beta\nabla_hu^{(n)}_s(\varphi),\dif L^{(n)}_s\>_0\right|^p\right)^{1/p}\\
&\leq C\|\Phi\|_\infty t^{-\frac{2}{\alpha}}
\left(\int^t_0\mE\|Q^{-1}_\beta\nabla_hu^{(n)}_s(\varphi)\|^\alpha_0\dif s\right)^{1/\alpha}\\
&\stackrel{(\ref{EE6})}{\leq} C\|\Phi\|_\infty t^{-\frac{2}{\alpha}}
\left(\int^t_0\mE\|A^{\frac{\theta}{2}}\nabla_hu^{(n)}_s(\varphi)\|^\alpha_0\dif s\right)^{1/\alpha}\\
&\stackrel{(\ref{EE7})}{\leq} C\|\Phi\|_\infty t^{-\frac{2}{\alpha}}
\left(\int^t_0s^{\frac{(\sigma-\theta)\alpha}{2}}\dif s\right)^{1/\alpha}\|h\|_\sigma\\
&\leq C\|\Phi\|_\infty t^{-\frac{1}{\alpha}-\frac{\theta-\sigma}{2}}\|h\|_\sigma,\ \ h\in\mH_n.
\end{align*}
From this, we in particular have
\begin{align}
|\mE \Phi(u^{(n)}_t(\varphi_1))-\mE \Phi(u^{(n)}_t(\varphi_2))|\leq C\|\Phi\|_\infty
t^{-\frac{1}{\alpha}-\frac{\theta-\sigma}{2}}\|\varphi_1-\varphi_2\|_\sigma,\ \ \varphi_1,\varphi_2\in\mH_n,\label{Es6}
\end{align}
where $C$ is independent of $n$.

(Step 3). In this step we prove that for any fixed $t>0$ and $\varphi\in\mH^0$,
\begin{align}
\lim_{n\to\infty}\|u^{(n)}_t(\Pi_n \varphi)-u_t(\varphi)\|_0=0,\ \ P-a.s.\label{ET9}
\end{align}
Set
$$
Z^{(n)}_t:=\int^t_0\mathrm{e}^{-(t-s)A}Q_\beta\dif L^{(n)}_s,\ \ w^{(n)}_t:=u^{(n)}_t-Z^{(n)}_t.
$$
Then
$$
w^{(n)}_t-w_t=\mathrm{e}^{-tA}(\Pi_n \varphi-\varphi)+\int^t_0\mathrm{e}^{-(t-s)A}(F_n(w^{(n)}_s+Z^{(n)}_s)-F(w_s+Z_s))\dif s,
$$
and
$$
\|w^{(n)}_t-w_t\|_\gamma\leq Ct^{-\frac{\gamma}{2}}\|\Pi_n \varphi-\varphi\|_0
+C\int^t_0(t-s)^{-\frac{\gamma+\gamma'}{2}}\|F_n(w^{(n)}_s+Z^{(n)}_s)-F(w_s+Z_s)\|_{-\gamma'}\dif s.
$$
Notice that by the definition of $F_n$,
\begin{align*}
\|F_n(w^{(n)}_s+Z^{(n)}_s)-F(w_s+Z_s)\|_{-\gamma'}
&\leq \|F\|_{\mathrm{Lip}}\Big(\|w^{(n)}_s-w_s\|_\gamma+\|(\Pi_n-I)Z_s\|_\gamma+\tfrac{1}{n}\Big)\no\\
&\quad+\|(\Pi_n-I)F(w_s+Z_s)\|_{-\gamma'}
\end{align*}
and
$$
\lim_{n\to\infty}\|(\Pi_n-I)Z_s\|_\gamma=0,\ \ \lim_{n\to\infty}\|(\Pi_n-I)F(w_s+Z_s)\|_{-\gamma'}=0.
$$
Since $F$ is bounded, by Fatou's lemma, we obtain
\begin{align}
\varlimsup_{n\to\infty}\|w^{(n)}_t-w_t\|_\gamma\leq C\int^t_0(t-s)^{-\frac{\gamma+\gamma'}{2}}
\varlimsup_{n\to\infty}\|w^{(n)}_s-w_s\|_\gamma\dif s,\label{EW2}
\end{align}
which then gives
\begin{align}
\varlimsup_{n\to\infty}\|w^{(n)}_t-w_t\|_\gamma=0\label{EW3}
\end{align}
as well as (\ref{ET9}).

(Step 4). For proving (\ref{For5}), we first assume $\Phi$ is continuous. In this case,
by taking limits for (\ref{Es6}), we obtain (\ref{For5}). For general bounded measurable $\Phi$, it follows by
a standard approximation.
\end{proof}

\section{Exponential Ergodicity of stochastic Burgers equations driven by $\alpha$-stable noises}

We first recall the following abstract form of Harris' theorem (cf. \cite[Theorem 4.2]{Ha-Ma-Sh}).
\bt\label{Ha}(Harris)
Let $\cP_t$ be a Markov semigroup over a Polish space $\mX$. We assume that for some Lyapunov function $V:\mX\to\mR_+$,
\begin{enumerate}[(i)]
\item there exist constants $C_V,\gamma,K_V>0$ such that for every $x\in\mX$ and $t>0$,
$$
\cP_tV(x)\leq C_V\mathrm{e}^{-\gamma t}V(x)+K_V;
$$
\item for every $R>0$, there exists a time $t>0$ and $\delta>0$ such that for all $x,y\in\mB^\mX_R$,
$$
\|\cP_t(x,\cdot)-\cP_t(y,\cdot)\|_{\mathrm{TV}}:=\sup_{\|\Phi\|\leq 1}|\cP_t\Phi(x)-\cP_t\Phi(y)|\leq 2-\delta,
$$
where $\|\cdot\|_{\mathrm{TV}}$ denotes the norm of total variation.
\end{enumerate}
Then $\cP_t$ has a unique invariant probability measure $\mu$ with
$$
\|\cP_t(x,\cdot)-\mu\|_{\mathrm{TV}}\leq C\mathrm{e}^{-\gamma_*t}(1+V(x))
$$
for some $C,\gamma_*>0$.
\et

In this section we shall use Theorems \ref{Th} and (\ref{Ha}) to prove
the exponential ergodicity of stochastic Burgers equations driven by $\alpha$-stable noises.
Let $\mH$ be the space of all square integrable functions on the torus $\mT=[0,2\pi)$
with vanishing mean values. Let $Au=-u''$ be the second order differential operator.
Then $A$ is a positive self-adjoint operator on $\mH$. Let $\lambda_{2k}:=\lambda_{2k+1}:=k^2$ and
$$
e_{2k}(x):=\pi^{-{1\over 2}}\cos(k x),\ \ e_{2k+1}(x):=\pi^{-{1\over 2}}\sin(k x).
$$
It is easy to see that $\{e_k, k\in\mN\}$ forms an orthogonal basis of $\mH$ and
$$
Ae_k=\lambda_k e_k,\  k\in\mN.
$$
In this case, let $\theta_0$ be defined by (\ref{ET2}), then
$$
\theta_0=\tfrac{1}{2}.
$$

Define a bilinear operator
$$
B(u,v):=uv',\ \ u,v\in\mH^1,
$$
and write
$$
B(u)=B(u,u).
$$
Consider the following stochastic Burgers equation driven by $L_t$:
\begin{align}
\dif u_t=[-Au_t-B(u_t)]\dif t+Q_\beta\dif L_t, \ u_0=\varphi\in \mH,\label{SBE}
\end{align}
where $Q_\beta$ denotes the intensity of the noise as above.

The main result of the paper is that
\bt\label{Main}
Let $\alpha\in(1,2)$. Assume that for some $\frac{3}{2}<\theta'\leq\theta<2$ and $\delta>0$,
\begin{align}
\delta k^{-\theta}\leq|\beta_k|\leq \delta^{-1}k^{-\theta'},\ \
\forall k\in\mN.\label{EW101}
\end{align}
\begin{enumerate}[(i)]
\item Let $Z_t:=\int^t_0 \mathrm{e}^{-(t-s)A}Q_\beta\dif L_s$.
Then for any $\varphi\in\mH$, there exists a unique $u_\cdot(\varphi)$ with
$$
u_\cdot-Z_\cdot\in C([0,\infty),\mH)\cap C((0,\infty),\mH^1)
$$
solving equation (\ref{SBE}). In particular, $(t,\varphi)\mapsto u_t(\varphi)$ is a Markov process on $\mH$. We write
$$
\cP_t\Phi(\varphi):=\mE\Phi(u_t(\varphi)).
$$
\item  $(\cP_t)_{t>0}$ is strong Feller, i.e.,
for any bounded measurable function $\Phi$ on $\mH$ and $t>0$, $\cP_t\Phi$ is a continuous function on $\mH$.

\item There exists a unique invariant probability measure $\mu$ on $\mH$ such that
\begin{align}
\|\cP_t(\varphi,\cdot)-\mu\|_{\mathrm{TV}}\leq C\mathrm{e}^{-\gamma_* t}(1+\|\varphi\|_0)\label{EQ1}
\end{align}
for some $C,\gamma_*>0$.
\end{enumerate}
\et
\begin{proof}
We divide the proof into four steps.

(Step 1). In view of $\theta'>\frac{3}{2}$, for any $\gamma\in(1,\theta'-\frac{1}{2})$,
by (\ref{Ep1}) we have
\begin{align}
\mE\left(\sup_{t\in[0,T]}\|Z_t\|_1\right)\leq
C_{\alpha}\left(\sum_kk^{2\gamma-2\theta'}\right) T^{\frac{1}{\alpha}}
\left(1+T^{\frac{\gamma-1}{2}}\right)<+\infty, \ \ T>0.\label{EW4}
\end{align}
Thus, (i) follows by Theorem \ref{Th3} below.\\

(Step 2). In this step, we prove the following claim: For given $R>0$,
there exist $T=T( R)\in(0,1]$ and $K_1=K_1(R,T)>0$, $K_2=K_2(\alpha,\theta')>0$ such that
for any bounded measurable function $\Phi$ on $\mH$, $\varphi_1,\varphi_2\in \mB^{\mH_1}_R$ and all $t\in(0,T]$,
\begin{align}
|\cP_t\Phi(\varphi_1)-\cP_t\Phi(\varphi_2)|\leq K_1
t^{-\frac{1}{\alpha}-\frac{\theta-1}{2}}\|\Phi\|_\infty\|\varphi_1-\varphi_2\|_1
+\frac{K_2t^{\frac{1}{\alpha}}}{R}.\label{EW9}
\end{align}

Consider the following truncated equation:
$$
\dif u^R_t=[-Au^R_t-B_R(u^R_t)]\dif t+Q_\beta\dif L_t, \ u^R_0=\varphi\in \mH^1,
$$
where
$$
B_R(u):=B(u)\cdot\chi(\|u\|_1/(5R)),
$$
and $\chi\in C^\infty(\mR, [0,1])$ satisfies
$$
\chi(r)=1,\ \forall |r|\leq 1;\ \ \chi(r)=0,\ \forall|r|>2.
$$
Define the stopping time
$$
\tau^R_\varphi(\omega):=\inf\Big\{t>0: \|u_t(\varphi;\omega)\|_1\geq 5R\Big\},\ \ \varphi\in \mB_{R}^{\mH^1},
$$
and let
$$
w_t(\varphi;\omega)=u_t(\varphi;\omega)-Z_t(\omega).
$$
Then we have
\begin{align}
\mP(\tau^R_{\varphi}\leq t)&=\mP\left(\sup_{s\in[0,t]}\|u_s(\varphi)\|_1\geq 5R\right)
\leq\mP\left(\sup_{s\in[0,t]}\|w_s(\varphi)\|_1+\sup_{s\in[0,t]}\|Z_s\|_1\geq 5R\right)\no\\
&\leq\mP\left(\sup_{s\in[0,t]}\|w_s(\varphi)\|_1\geq 4R, \sup_{s\in[0,t]}\|Z_s\|_1\leq R\right)
+\mP\left(\sup_{s\in[0,t]}\|Z_s\|_1>R\right).\label{EU1}
\end{align}
By Theorem \ref{Th3} below, there exists a time $T=T(R)\in(0,1]$
such that for any $\varphi\in \mB_R^{\mH^1}$ and $t\in(0,T]$,
if $\sup_{s\in[0,t]}\|Z_s(\omega)\|_1\leq R$, then
\begin{align}
\sup_{s\in[0,t]}\|w_s(\varphi;\omega)\|_1\leq 3R.\label{ER3}
\end{align}
Hence, by (\ref{EU1}) and Chebychev's inequality, we have for any $\varphi\in \mB_R^{\mH^1}$,
\begin{align}
\mP(\tau^R_{\varphi}\leq t)\leq \mP\left(\sup_{s\in[0,t]}\|Z_s\|_1>R\right)
\leq\frac{\mE\left(\sup_{s\in[0,t]}\|Z_s\|_1\right)}{R}
\stackrel{(\ref{EW4})}{\leq}\frac{C_{\alpha,\theta'}t^{\frac{1}{\alpha}}}{R}.\label{EU2}
\end{align}
On the other hand, by the uniqueness of solutions, we have
$$
u_t(\varphi)=u^R_t(\varphi), \forall t\in[0,\tau^R_{\varphi}).
$$
Thus, if we choose $\sigma=\gamma=1$ and $\gamma'=0$ in Theorem \ref{Th}, then by (\ref{For5}),
we have for any $t\in(0,T]$ and $\varphi_1,\varphi_2\in \mB_R^{\mH^1}$,
\begin{align*}
|\cP_t\Phi(\varphi_1)-\cP_t\Phi(\varphi_2)|&\leq\Big|\mE\Big(\Phi(u_t(\varphi_1)); \tau^R_{\varphi_1}>t\Big)-
\mE\Big(\Phi(u_t(\varphi_2)); \tau^R_{\varphi_2}>t\Big)\Big|
+\mP(\tau^R_{\varphi_1}\leq t)+\mP(\tau^R_{\varphi_2}\leq t)\\
&=\Big|\mE\Big(\Phi(u^R_t(\varphi_1)); \tau^R_{\varphi_1}>t\Big)-
\mE\Big(\Phi(u^R_t(\varphi_2)); \tau^R_{\varphi_2}>t\Big)\Big|
+\mP(\tau^R_{\varphi_1}\leq t)+\mP(\tau^R_{\varphi_2}\leq t)\\
&\leq\Big|\mE\Big(\Phi(u^R_t(\varphi_1))\Big)-\mE\Big(\Phi(u^R_t(\varphi_2))\Big)\Big|
+2\mP(\tau^R_{\varphi_1}\leq t)+2\mP(\tau^R_{\varphi_2}\leq t)\\
&\leq K_1t^{-\frac{1}{\alpha}-\frac{\theta-1}{2}}\|\Phi\|_\infty\|\varphi_1-\varphi_2\|_1
+2\mP(\tau^R_{\varphi_1}\leq t)+2\mP(\tau^R_{\varphi_2}\leq t),
\end{align*}
which together with (\ref{EU2}) gives (\ref{EW9}).\\

(Step 3). In this step, we prove (ii). Let $\Phi$ be a bounded measurable function
on $\mH$. Let us first show that for any $t>0$,
\begin{align}
\varphi\mapsto\cP_t\Phi(\varphi)\mbox{ is continuous on $\mH^1$}.\label{ER7}
\end{align}
Let $\{\varphi_n\}\subset\mB^{\mH^1}_R$ converge to $\varphi$ in $\mH^1$.
Let $T=T(R)\in(0,1]$ be as in Step 2. For fixed $t>0$, by (\ref{EW9}) we have
\begin{align*}
|\cP_t\Phi(\varphi_n)-\cP_t\Phi(\varphi)|
&=|\cP_{t\wedge T}\cP_{t-t\wedge T}\Phi(\varphi_n)-\cP_{t\wedge T}\cP_{t-t\wedge T}\Phi(\varphi)|\\
&\leq K_1 (t\wedge T)^{-\frac{1}{\alpha}-\frac{\theta-1}{2}}\|\cP_{t-t\wedge T}\Phi\|_\infty\|\varphi_n-\varphi\|_1
+\frac{K_2(t\wedge T)^{\frac{1}{\alpha}}}{R}\\
&\leq K_1 (t\wedge T)^{-\frac{1}{\alpha}-\frac{\theta-1}{2}}\|\Phi\|_\infty\|\varphi_n-\varphi\|_1
+\frac{K_2}{R},
\end{align*}
where $K_1=K_1(R,T)$ and $K_2=K_2(\alpha,\theta')$.
First letting $n\to\infty$ and then $R\to\infty$, we obtain
$$
\lim_{n\to\infty}|\cP_t\Phi(\varphi_n)-\cP_t\Phi(\varphi)|=0.
$$

Next we prove that
\begin{align}
\varphi\mapsto\cP_t\Phi(\varphi)\mbox{ is continuous on $\mH$}.\label{ER1}
\end{align}
For $R>0$, define
$$
\Omega_R:=\left\{\omega: \sup_{s\in[0,1]}\|Z_s(\omega)\|_1\leq R\right\}.
$$
By Theorem \ref{Th3} again, there exists a time $T=T(R)\in(0,1)$ such that
for any $\omega\in\Omega_R$ and all $\varphi_1,\varphi_2\in \mB_R^\mH$ and $t\in(0,T]$,
\begin{align}
\|u_t(\varphi_1;\omega)-u_t(\varphi_2;\omega)\|_1=\|w_t(\varphi_1;\omega)-w_t(\varphi_2;\omega)\|_1
\leq 2t^{-\frac{1}{2}}\|\varphi_1-\varphi_2\|_0.\label{ER4}
\end{align}
Let $\{\varphi_n\}\subset\mB^\mH_R$ converge to $\varphi$ in $\mH$.
For any $t>0$, we have
\begin{align*}
|\cP_t\Phi(\varphi_n)-\cP_t\Phi(\varphi)|&
=\left|\mE\Big((\cP_{t-t\wedge T}\Phi)(u_{t\wedge T}(\varphi_n))\Big)
-\mE\Big((\cP_{t-t\wedge T}\Phi)(u_{t\wedge T}(\varphi))\Big)\right|\\
&\leq \left|\mE\Big((P_{t-t\wedge T}\Phi)(u_{t\wedge T}(\varphi_n))
-(P_{t-t\wedge T}\Phi)(u_{t\wedge T}(\varphi)); \Omega_R\Big)\right|+2\mP(\Omega_R^c),
\end{align*}
which together with  (\ref{ER7}) and (\ref{ER4}) yields (\ref{ER1}) by first letting $n\to\infty$ and then $R\to\infty$.\\

(Step 4). In this step, we prove (iii). Take $V(\varphi)=1+\|\varphi\|_0$.
Let us first check (i) of Theorem \ref{Ha}.
Arguing as deriving (1.2) of \cite{Do-Xu-Zh} and taking $\theta=1$ therein, we have
\begin{align*}
\mE\left(\sup_{s\in[0,t]}\|u_s\|_0\right)+
\mE\left(\int^t_0\frac{\|u_s\|^2_1}{(\|u_s\|^2_0+1)^{1/2}}\dif s\right)
\leq C(1+\|\varphi\|_0+t). 
\end{align*}
which, together with the spectral gap inequality $\|u\|_0 \le \|u\|_1$, implies
\begin{align*}
\mE\left(\|u_t\|_0+1\right)+
\mE\left(\int^t_0\frac{\|u_s\|^2_0+1}{(\|u_s\|^2_0+1)^{1/2}}\dif s\right)
\leq C(1+\|\varphi\|_0+t). 
\end{align*}
From this, we get
\begin{align*}
\mE V(u_t)\leq -\frac{1}{2}
\int^t_0 \mE V(u_s)\dif s+CV(\varphi)+Ct,
\end{align*}
which implies that
\begin{equation}
\mE V(u_t) \le C \mathrm{e}^{-\frac 12 t}V(\varphi)+2C, \ \ \ \forall \ t>0.
\end{equation}

Next we check (ii) of  Theorem \ref{Ha}. Fix $R>0$. Let $\eps, t_0>0$, to be determined later. Define
$$
\Omega^\eps_{t_0}:=\left\{\omega: \sup_{s\in[0,t_0+1]}\|Z_s(\omega)\|_1\leq\eps\right\}.
$$
By (\ref{ER9}) below, one can choose $\eps_0:=\frac{1}{2C_1}\wedge\frac{1}{\sqrt{2C_2+1}}$
and $t_0\geq2\log(R^2/\eps^4)$ so that for each $\eps\in(0,\eps_0]$,
all $\omega\in\Omega^\eps_{t_0}$, $\varphi\in \mB^{\mH}_R$ and $t\in[t_0,t_0+1]$,
\begin{align*}
\|w_t(\varphi,\omega)\|_0^2&\leq \|\varphi\|_0^2\mathrm{e}^{(C_1\eps-1)t}
+C_2\eps^4\int^t_0\mathrm{e}^{(C_1\eps-1)(t-s)}\dif s\\
&\leq R^2\mathrm{e}^{-t/2}+2C_2\eps^4\leq (2C_2+1)\eps^4\leq\eps^2.
\end{align*}
By using Theorem \ref{Th3} again with $R=\eps$ and starting from $t_0$ therein,
there exits a time $t_1\in(t_0,t_0+1]$ such that
for all $t\in(t_0,t_1]$ and $\eps\in(0,\eps_0]$,
$$
\|w_t(\varphi,\omega)\|_1\leq (t-t_0)^{-\frac{1}{2}}(2\|w_{t_0}(\varphi,\omega)\|_0)\leq(t-t_0)^{-\frac{1}{2}}(2\eps).
$$
In particular, for each $\eps\in(0,\eps_0)$, all $\omega\in\Omega^\eps_{t_0}$ and $\varphi\in \mB^{\mH}_R$,
$$
\|w_{t_1}(\varphi,\omega)\|_1\leq 2(t_1-t_0)^{-\frac{1}{2}}\eps.
$$
For $\varphi_1,\varphi_2\in \mB^{\mH}_R$, if we define
$$
A^\eps_{t_1}(\varphi_1,\varphi_2):=\Big\{\omega: \|w_{t_1}(\varphi_1,\omega)\|_1
+\|w_{t_1}(\varphi_2,\omega)\|_1\leq4(t_1-t_0)^{-\frac{1}{2}}\eps\Big\},
$$
then from the above implication, one has
\begin{align}
\Omega^\eps_{t_0}\subset A^\eps_{t_1}(\varphi_1,\varphi_2).\label{EU3}
\end{align}
Now by definition, for any $t_2\in(t_1,t_0+1)$ with $t_2-t_1$ being small, we have
\begin{align*}
\|\cP_{t_2}(\varphi_1,\cdot)-\cP_{t_2}(\varphi_2,\cdot)\|_{\mathrm{TV}}
&:=\sup_{\|\Phi\|_\infty\leq 1}|\cP_{t_2}\Phi(\varphi_1)-\cP_{t_2}\Phi(\varphi_2)|\\
&=\sup_{\|\Phi\|_\infty\leq 1}|\cP_{t_1}\cP_{t_2-t_1}\Phi(\varphi_1)-\cP_{t_1}\cP_{t_2-t_1}\Phi(\varphi_2)|\\
&=\sup_{\|\Phi\|_\infty\leq 1}\Big|\mE\Big(\cP_{t_2-t_1}\Phi(u_{t_1}(\varphi_1))
-\cP_{t_2-t_1}\Phi(u_{t_1}(\varphi_2))\Big)\Big|\\
&\leq \sup_{\|\Phi\|_\infty\leq 1}\Big|\mE\Big(\cP_{t_2-t_1}\Phi(u_{t_1}(\varphi_1))
-\cP_{t_2-t_1}\Phi(u_{t_1}(\varphi_2)); A^\eps_{t_1}(\varphi_1,\varphi_2)\Big)\Big|\\
&\qquad+2\Big(1-\mP(A^\eps_{t_1}(\varphi_1,\varphi_2))\Big).
\end{align*}
Noticing that on $A^\eps_{t_1}(\varphi_1,\varphi_2)$,
$$
\|u_{t_1}(\varphi_1)-u_{t_1}(\varphi_2)\|_1=
\|w_{t_1}(\varphi_1)-w_{t_1}(\varphi_2)\|_1\leq4(t_1-t_0)^{-\frac{1}{2}}\eps,
$$
by (\ref{EW9}), we further have for all $\eps\in(0,\eps_0)$,
\begin{align*}
\|\cP_{t_2}(\varphi_1,\cdot)-\cP_{t_2}(\varphi_2,\cdot)\|_{\mathrm{TV}}
&\leq \Bigg(4K_1 (t_2-t_1)^{-\frac{1}{\alpha}-\frac{\theta-1}{2}}(t_1-t_0)^{-\frac{1}{2}}\eps
+\frac{K_2(t_2-t_1)^{\frac{1}{\alpha}}}{R}\Bigg)\\
&\qquad\times\mP(A^\eps_{t_1}(\varphi_1,\varphi_2)+2\Big(1-\mP(A^\eps_{t_1}(\varphi_1,\varphi_2))\Big)\\
&=2-\Bigg(2-4K_1(t_2-t_1)^{-\frac{1}{\alpha}-\frac{\theta-1}{2}}(t_1-t_0)^{-\frac{1}{2}}\eps
-\frac{K_2(t_2-t_1)^{\frac{1}{\alpha}}}{R}\Bigg)\\
&\qquad\qquad\times\mP(A^\eps_{t_1}(\varphi_1,\varphi_2)).
\end{align*}
Choosing first $t_2\in(t_1,t_0+1)$ so that
$$
\frac{K_2(t_2-t_1)^{\frac{1}{\alpha}}}{R}\leq \frac{1}{2},
$$
and then $\eps\in(0,\eps_0)$ so that
$$
4K_1(t_2-t_1)^{-\frac{1}{\alpha}-\frac{\theta-1}{2}}(t_1-t_0)^{-\frac{1}{2}}\eps\leq\tfrac{1}{2},
$$
we finally obtain that for all $\varphi_1,\varphi_2\in\mB^\mH_R$,
$$
\|\cP_{t_2}(\varphi_1,\cdot)-\cP_{t_2}(\varphi_2,\cdot)\|_{\mathrm{TV}}\leq 2-
\mP(A^\eps_{t_1}(\varphi_1,\varphi_2))\stackrel{(\ref{EU3})}{\leq} 2-\mP(\Omega^\eps_{t_0}).
$$
The condition (ii) of Theorem \ref{Ha} is thus verified by (\ref{Ep11}), and
(iii) follows by Theorem \ref{Ha}. The whole proof is complete.
\end{proof}

\section{Appendix: A study of deterministic Burgers equation}

In this appendix we study the following deterministic Burgers equation:
\begin{align}
\dot w_t=-Aw_t-B(w_t+Z_t), \ w_0=\varphi\in \mH^0,\label{SNSE}
\end{align}
where $t\mapsto Z_t$ is a bounded measurable function on $\mH^1$.

Recall the following estimate about the bilinear form $B(u,v)$ (see \cite[Lemma 2.1]{Te}):
\begin{align}
\<B(u,v),w\>_0\leq C\|u\|_{\sigma_1}\|v\|_{\sigma_2+1}\|w\|_{\sigma_3},\ \ \sigma_1+\sigma_2+\sigma_3>1/2,\label{EE2}
\end{align}
where $C$ only depends on $\sigma_1,\sigma_2,\sigma_3$.
Let $\mM_T$ be the Banach space defined by
$$
\mM_T:=\left\{u\in C([0,T];\mH)\cap C((0,T];\mH^1):
\|u\|_{\mM_T}:=\sup_{t\in[0,T]}(\|u_t\|_0\vee (t^{\frac{1}{2}}\|u_t\|_1))<+\infty\right\}.
$$
We have
\bt\label{Th3}
For given $R>0$, there exists a time $T=T(R)\in(0,1]$, which is increasing
as $R\downarrow 0$,  such that if $\sup_{t\in[0,T]}\|Z_t\|_1\leq R$, then
\begin{enumerate}[(i)]
\item for any $\varphi\in \mB_{R}^\mH$,
there is a unique $w=w(\varphi)\in \mB_{2R}^{\mM_T}$ satisfying that for all $t\in[0,T]$,
\begin{align}
w_t=\mathrm{e}^{-tA}\varphi-\int^t_0 \mathrm{e}^{-(t-s)A}B(w_s+Z_s)\dif s;\label{EE5}
\end{align}
\item for any $\varphi_1,\varphi_2\in \mB_{R}^\mH$,
\begin{align}
\|w_\cdot(\varphi_1)-w_\cdot(\varphi_2)\|_{\mM_T}\leq 2\|\varphi_1-\varphi_2\|_0;\label{EE3}
\end{align}
\item for any $\varphi\in \mB_{R}^{\mH^1}$ and $t\in[0,T]$,
\begin{align}
\|w_t(\varphi)\|_1\leq 3R.\label{EE4}
\end{align}
\end{enumerate}
Moreover, there are two constants $C_1,C_2>0$ such that for any $\varphi\in\mH$ and all $t\geq 0$,
\begin{align}
\|w_t\|_0^2\leq \|\varphi\|_0^2\mathrm{e}^{\int^t_0(C_1\|Z_s\|^2_1-1)\dif s}
+C_2\int^t_0\mathrm{e}^{\int^t_s(C_1\|Z_r\|^2_1-1)\dif r}\|Z_s\|_1^4\dif s.\label{ER9}
\end{align}
In particular, for any $\varphi\in\mH$, there exists a unique $w_\cdot(\varphi)\in \cup_{T>0}\mM_T$ satisfying (\ref{EE5}).
\et
\begin{proof}
We use the fixed point argument. Fix $\varphi\in \mB_{R}^\mH$. Define a nonlinear map on $\mM_T$ by
$$
\cM(w)_t:=\mathrm{e}^{-tA}\varphi-\int^t_0 \mathrm{e}^{-(t-s)A}B(w_s+Z_s)\dif s.
$$
We want to show that for some $T:=T(R)\leq 1$,
$$
\cM\mbox{ is a contraction operator on $\mB_{2R}^{\mM_T}$}.
$$
Fix $\sigma\in(\frac{1}{2},1)$. For $w\in\mB_{2R}^{\mM_T}$, by (\ref{Semi}) and (\ref{EE2}), we have
for all $t\leq T$,
\begin{align*}
\|\cM (w)_t\|_0&\leq\|\varphi\|_0+C_\sigma\int^t_0(t-s)^{-\frac{\sigma}{2}}\|B(w_s+Z_s)\|_{-\sigma}\dif s\\
&\leq R+C_\sigma\int^t_0(t-s)^{-\frac{\sigma}{2}}(\|w_s\|_0+\|Z_s\|_0)(\|w_s\|_1+\|Z_s\|_1)\dif s\\
&\leq R+C_\sigma R\int^t_0(t-s)^{-\frac{\sigma}{2}}(s^{-\frac{1}{2}} R+ R)\dif s\\
&\leq R+C_\sigma R\Big(t^{\frac{1-\sigma}{2}} R+t^{1-\frac{\sigma}{2}} R\Big)
\leq R+C_\sigma R^2t^{\frac{1-\sigma}{2}},
\end{align*}
where $C_\sigma$ only depends on $\sigma$. Similarly, we also have
\begin{align*}
\|\cM(w)_t\|_1&\leq t^{-\frac{1}{2}}\|\varphi\|_0+C_\sigma\int^t_0(t-s)^{-\frac{1+\sigma}{2}}
\|B(w_s+Z_s)\|_{-\sigma}\dif s\\
&\leq t^{-\frac{1}{2}}R+C_\sigma R^2\int^t_0(t-s)^{-\frac{1+\sigma}{2}}s^{-\frac{1}{2}}\dif s\\
&\leq t^{-\frac{1}{2}} R+C_\sigma R^2t^{-\frac{\sigma}{2}}.
\end{align*}
Hence,
$$
\|\cM(w)\|_{\mM_T}\leq  R+C_\sigma R^2T^{\frac{1-\sigma}{2}}.
$$
If we choose $T\leq(C_\sigma R)^{-\frac{2}{1-\sigma}}\wedge 1=:T_1$, then $\|\cM(w)\|_{\mM_T}\leq 2 R$, and
$$
\mbox{$\cM$ maps $\mB^{\mM_T}_{2R}$ into $\mB^{\mM_T}_{2R}$.}
$$
On the other hand, for $w,v\in\mB_{2R}^{\mM_T}$,  by (\ref{EE2}) again, we have
\begin{align*}
\|\cM(w)_t-\cM(v)_t\|_0&\leq C_\sigma\int^t_0(t-s)^{-\frac{\sigma}{2}}\|B(w_s+Z_s)
-B(v_s+Z_s)\|_{-\sigma}\dif s\\
&\leq C_\sigma\int^t_0(t-s)^{-\frac{\sigma}{2}}\|w_s-v_s\|_0(\|w_s\|_1+\|Z_s\|_1)\dif s\\
&+C_\sigma\int^t_0(t-s)^{-\frac{\sigma}{2}}\|w_s-v_s\|_1(\|v_s\|_0+\|Z_s\|_0)\dif s\\
&\leq C_\sigma\sup_{s\in[0,t]}\|w_s-v_s\|_0\int^t_0(t-s)^{-\frac{\sigma}{2}}(s^{-\frac{1}{2}} R+ R)\dif s\\
&+C_\sigma\sup_{s\in[0,t]}s^{\frac{1}{2}}\|w_s-v_s\|_1\int^t_0(t-s)^{-\frac{\sigma}{2}}s^{-\frac{1}{2}} R\dif s\\
&\leq C_\sigma R\|w-v\|_{\mM_t}t^{\frac{1-\sigma}{2}},
\end{align*}
and
$$
\|\cM (w)_t-\cM (v)_t\|_1\leq C_\sigma R\|w-v\|_{\mM_t}t^{-\frac{\sigma}{2}}.
$$
Hence,
$$
\|\cM (w)-\cM(v)\|_{\mM_T}\leq C_\sigma R\|w-v\|_{\mM_T}T^{\frac{1-\sigma}{2}}.
$$
Letting $T\leq\frac{1}{2C_\sigma R}\wedge T_1=:T_2$, we obtain
\begin{align}
\|\cM (w)-\cM (v)\|_{\mM_T}\leq \tfrac{1}{2}\|w-v\|_{\mM_T}.\label{EU4}
\end{align}
The existence and uniqueness for equation (\ref{EE5}) follow by the fixed point theorem.
Moreover, as in estimating (\ref{EU4}), we also have (\ref{EE3}).

Next we prove (\ref{EE4}). As above, by (\ref{Semi}) and (\ref{EE2}), we have
\begin{align*}
\|w_t\|_1&\leq\|\mathrm{e}^{-tA}\varphi\|_1+\int^t_0 \|\mathrm{e}^{-(t-s)A}B(w_s+Z_s)\|_1\dif s\\
&\leq\|\varphi\|_1+C_\sigma\int^t_0 (t-s)^{-\frac{1+\sigma}{2}}\|B(w_s+Z_s)\|_{-\sigma}\dif s\\
&\leq\|\varphi\|_1+C_\sigma\int^t_0 (t-s)^{-\frac{1+\sigma}{2}}(\|w_s\|_0+\|Z_s\|_0)(\|w_s\|_1+\|Z_s\|_1)\dif s\\
&\leq\|\varphi\|_1+C_\sigma R\int^t_0 (t-s)^{-\frac{1+\sigma}{2}}(\|w_s\|_1+ R)\dif s\\
&\leq\|\varphi\|_1+C_\sigma R t^{\frac{1-\sigma}{2}}\left(\sup_{s\in[0,t]}\|w_s\|_1+ R\right).
\end{align*}
From this, one sees that for $t\leq(2(C_\sigma R))^{-\frac{2}{1-\sigma}}\wedge T_2$,
$$
\sup_{s\in[0,t]}\|w_s\|_1\leq 2\|\varphi\|_1+2C_\sigma R^2 t^{\frac{1-\sigma}{2}}\leq 3 R.
$$

We now prove (\ref{ER9}). Notice that
$$
\p_t\|w_t\|_0^2=-2\|w_t\|_1^2+2\<B(w_t+Z_t),w_t\>_0.
$$
Since $\<B(w,w),w\>_0=0$, we have
$$
\<B(w_t+Z_t),w_t\>_0=\<B(w_t,Z_t),w_t\>_0+\<B(Z_t,Z_t),w_t\>_0+\<B(Z_t,w_t),w_t\>_0.
$$
Thus, by (\ref{EE2}) and Young's inequality, we have
\begin{align*}
2|\<B(w_t+Z_t),w_t\>_0|&\leq C\|w_t\|_1\|Z_t\|_1\|w_t\|_0+C\|Z_t\|^2_1\|w_t\|_1\\
&\leq\|w_t\|^2_1+C_1\|w_t\|_0^2\|Z_t\|^2_1+C_2\|Z_t\|_1^4.
\end{align*}
Hence, by $\|w\|_0\leq \|w\|_1$, we obtain
$$
\p_t\|w_t\|_0^2\leq (C_1\|Z_t\|^2_1-1)\|w_t\|_0^2+C_2\|Z_t\|_1^4,
$$
which implies (\ref{ER9}) by solving this differential inequality.
\end{proof}

{\bf Acknowledgements:}

We would like to thank Vahagn Nersesyan for stimulating discussions on exponential mixing.
Zhao Dong is supported by 973 Program (2011CB808000), Science Fund for Creative Research Groups (10721101),
NSFs of China (Nos. 10928103, 11071008). Xicheng Zhang is supported by NSFs of China (No. 11271294)
and Program for New Century Excellent Talents in University.

\end{document}